\documentclass[11pt,oneside]{amsart}
\usepackage[english]{babel}

\makeatletter
\def\dom{\mathop{\operatorname{dom}}}
\def\vrai#1{\mathop{\operatorname{vrai\:#1}}}
\def\supp{\mathop{\operatorname{supp}}}
\def\EMPTY@@@{}
\def\subsection#1{\refstepcounter{subsection}\par%
	\vskip 0.3cm\noindent{\bfseries\upshape\thesubsection .}%
	~{\bfseries\upshape #1}}
\def\subsubsection#1{\refstepcounter{subsubsection}\par%
	\vskip 0.3cm\noindent{\bfseries\upshape\thesubsubsection .}%
	~{\bfseries\upshape #1}}
\def\Wo{{\mathpalette\Wo@{}}W}
\def\Wo@#1{\setbox0\hbox{$#1 W$}\dimen@\ht0\dimen@ii\wd0\raise0.65\dimen@%
\rlap{\kern0.35\dimen@ii$#1{}^\circ$}}
\makeatother

\title[The extrema of Sturm--Liouville eigenvalue]{The Sturm--Liouville problem
with singular potential and the extrema of the first eigenvalue}
\author{E.\,S.~Karulina}
\thanks{The first author is supported by the Russian Foundation for Basic
Researches, grant No~11-01-00989.}
\address{Moscow State University of Economics, Statistics and Informatics}
\email{KarulinaES@yandex.ru}
\author{A.\,A.~Vladimirov}
\thanks{The second author is supported by the Russian
Foundation for Basic Researches, grant No~10-01-00423.}
\address{Dorodnitsyn Computing Center of RAS}
\email{vladimi@mech.math.msu.su}

\begin{document}
\begin{abstract}
We get the infima and suprema of the first eigenvalue of the problem
\begin{gather*}
	-y''+qy=\lambda y,\\
	\left\{\begin{aligned}
		y'(0)-k_0^2y(0)=0,\\ y'(1)+k_1^2y(1)=0,
	\end{aligned}\right.
\end{gather*}
where \(q\) belongs to the set of constant-sign summable functions on \([0,1]\)
such that
\[
	\int_0^1 q\,dx=1 \text{ or }\int_0^1 q\,dx=-1.
\]
\end{abstract}
\maketitle

\section{Introduction}\label{par:1}
\subsection{}
Consider the Sturm--Liouville problem
\begin{gather}\label{eq:eq1}
	-y''+(q-\lambda)y=0,\\ \label{eq:eq2}
	\left\{\begin{aligned}
		y'(0)-k_0^2y(0)=0,\\ y'(1)+k_1^2y(1)=0,
	\end{aligned}\right.
\end{gather}
where real coefficients \(k_0\geq 0\) and \(k_1\geq k_0\) are fixed, the solution \(y\)
belongs to space \(W_1^2[0,1]\), equality \eqref{eq:eq1} is considered as holding
almost everywhere at \([0,1]\), and the potential \(q\in L_1[0,1]\)
is a constant-sign function such that one of the integral conditions holds:
\begin{equation}\label{eq:int}
	\int_0^1 q\,dx=1 \;\text{ or } \int_0^1 q\,dx=-1.
\end{equation}
The aim of this paper is to get the infima and suprema of the first eigenvalue
of problem \eqref{eq:eq1}--\eqref{eq:int}.

\subsection{}\label{pt:vrai1.2}
Problem \eqref{eq:eq1}--\eqref{eq:int} is a partial case of problem
\eqref{eq:eq1}, \eqref{eq:eq2} with \(q\in A_{\gamma}\) or \(-q\in A_{\gamma}\),
where \(\gamma\in\mathbb R\setminus\{0\}\) and
\begin{equation}\label{eq:Agamma}
	A_{\gamma}\rightleftharpoons\left\{q\in L_1[0,1]\::\: q(x)\geq 0
		\text{\ a.\,e.,\ and\ }\int_0^1 q^{\gamma}\,dx=1\right\}.
\end{equation}
Denote by \(\lambda_1(q)\) the minimal eigenvalue of problem \eqref{eq:eq1} or
\begin{equation}\label{eq:EK}
	-y''-\lambda qy=0
\end{equation}
with some self-adjoint boundary conditions. Consider for each \(\gamma\in
\mathbb R\setminus\{0\}\) four values \(m_{\gamma}^{\pm}\rightleftharpoons
\inf_{q\in A_{\gamma}}\lambda_1(\pm q)\) and \(M_{\gamma}^{\pm}
\rightleftharpoons\sup_{q\in A_{\gamma}}\lambda_1(\pm q)\). The estimates of
\(m_{\gamma}^+\) and \(M_{\gamma}^+\) for equation \eqref{eq:EK} with
Dirichlet boundary conditions were obtained in \cite{EK:1996}. The analogous
results about Dirichlet problem for equation \eqref{eq:eq1} were obtained
in \cite{VS:2003}, \cite{Ez:2005}. In \cite{Mur:2005} problem \eqref{eq:EK},
\eqref{eq:eq2} was studied.

The values \(m_{\gamma}^+\) and \(M_{\gamma}^+\) for problem \eqref{eq:eq1},
\eqref{eq:eq2} with \(q\in A_\gamma\) were considered by one of the authors
in \cite{Kar:2011} for all \(\gamma\neq 0\). The most detailed and precise results
were obtained for the case \(\gamma\ne 1\).

The case \(\gamma=1\) is in some kind special. In \cite{Ez:2005} and
\cite{Kar:2011}, for \eqref{eq:eq1} with various boundary conditions, the precise
results for \(M_1^+\) were obtained by the method quite different from used for
\(\gamma\ne 1\). In \cite{Kar:2011} for \(m_1^+\) only inequality \(m_1^+\geq
1/4\) was obtained. In \cite{Ez:2005} for \(m_1^-\) it was proved that this infimum
is attained at the non-summable potential \(q^*=-\boldsymbol\delta_{1/2}\).

In this paper we extend the class of considered potentials from \(L_1[0,1]\)
to space \(W_2^{-1}[0,1]\) (see \cite{SaSh:2003} and \ref{pt:1.2} later). Space
\(W_2^{-1}[0,1]\), in particular, contains a Dirac delta function
\(\boldsymbol\delta_{\zeta}\) with support located at an arbitrary point
\(\zeta\in [0,1]\). This generalization of the problem lets us to get the precise
description of \(M_1^{-}\) and \(m_1^{\pm}\) and to prove that they are attained
at the potentials from the extended class.

\subsection{}
The main results of the paper are following four theorems:

\subsubsection{}\label{prop:M1+}
{\textbf{Theorem.~}\itshape By definition, put
\begin{equation}\label{eq:ab}
	\alpha_{\mu}\rightleftharpoons\dfrac{1}{\sqrt{\mu}}\arctan
		\dfrac{k_0^2}{\sqrt{\mu}},
	\qquad\beta_{\mu}\rightleftharpoons\dfrac{1}{\sqrt{\mu}}\arctan
		\dfrac{k_1^2}{\sqrt{\mu}}.
\end{equation}
Then \(M^+_1\) is a unique solution to the equation
\begin{equation}\label{eq:M1}
	1-\alpha_{\mu}-\beta_{\mu}=\mu^{-1}
\end{equation}
and is attained at the potential \(q^*\in L_1[0,1]\) such that
\[
	\belowdisplayskip=-0.4cm
	q^*(x)=\left\{\begin{aligned}&M^+_1&
		&\text{for }x\in[\alpha_{M_1^+},1-\beta_{M_1^+}],\\
		&0&&\text{otherwise}.
	\end{aligned}\right.
\]
}

\subsubsection{}\label{prop:M1-}
{\textbf{Theorem.~}\itshape If \(k_0^2+k_1^2\leq 1\), then \(M^{-}_1=
k_0^2+k_1^2-1\) and is attained at the potential
\[
	q^*\rightleftharpoons -k_0^2\boldsymbol{\delta}_0-
		k_1^2\boldsymbol{\delta}_1-(1-k_0^2-k_1^2).
\]

If \(k_0^2+k_1^2\geq 1\) and \(k_1^2-k_0^2\leq 1\), then \(M^{-}_1\)
is the minimal eigenvalue of the problem
\begin{gather}\label{eq:propM1-1}
	-y''=\lambda y,\\ \label{eq:propM1-2}
	2y'(0)-(k_0^2+k_1^2-1)y(0)=2y'(1)+(k_0^2+k_1^2-1)y(1)=0
\end{gather}
and is attained at the potential
\[
	q^*\rightleftharpoons -(1+k_0^2-k_1^2)\boldsymbol{\delta}_0/2-
		(1-k_0^2+k_1^2)\boldsymbol{\delta}_1/2.
\]

If \(k_1^2-k_0^2\geq 1\), then \(M^{-}_1\) is the minimal eigenvalue
of problem \eqref{eq:propM1-1} with
\begin{equation}\label{eq:propM1-3}
	y'(0)-k_0^2y(0)=y'(1)+(k_1^2-1)y(1)=0
\end{equation}
and is attained at the potential \(q^*\rightleftharpoons-\boldsymbol{\delta}_1\).
}

\subsubsection{}\label{prop:m1+}
{\textbf{Theorem.~}\itshape \(m^+_1\) is the minimal eigenvalue of problem
\eqref{eq:propM1-1} with
\begin{equation}\label{eq:propm1+}
	y'(0)-k_0^2y(0)=y'(1)+(k_1^2+1)y(1)=0
\end{equation}
and is attained at the potential \(q^*\rightleftharpoons\boldsymbol{\delta}_1\).
}

\subsubsection{}\label{prop:m1-}
{\textbf{Theorem.~}\itshape If for some \(\mu\geq-k_0^4\) and some \(\zeta\in
(0,1)\) the problem
\begin{gather}\label{eq:grmm1}
	-y''=\mu y\quad\text{at }(0,\zeta)\cup(\zeta,1),\\ \label{eq:grmm2}
	y'(0)-k_0^2y(0)=2y'(\zeta-0)-y(\zeta)=2y'(\zeta+0)+y(\zeta)=
		y'(1)+k_1^2y(1)=0
\end{gather}
has a continuous positive solution, then \(m^{-}_1=\mu\) and \(m^{-}_1\)
is attained at the potential \(q^*\rightleftharpoons-\boldsymbol{\delta}_{\zeta}\).
Otherwise \(m^{-}_1\) is the minimal eigenvalue of problem \eqref{eq:propM1-1}
with
\[
	y'(0)-(k_0^2-1)y(0)=y'(1)+k_1^2y(1)=0
\]
and is attained at the potential \(q^*\rightleftharpoons -\boldsymbol{\delta}_0\).
}

\medskip
Some additional remarks on solvability of boundary problem \eqref{eq:grmm1},
\eqref{eq:grmm2} will be given in subsection \ref{pt:3.5}.

\subsection{}
Let us give some examples that illustrate the theorems from previous subsection.
In the case \(k_0=k_1=0\) we get \(m_1^+=\lambda_1(\boldsymbol\delta_1)=
0.740174(\pm 10^{-6})\). In the case \(k_0^2=k_1^2>1/2\) we get \(m_1^-=
\lambda_1(-\boldsymbol\delta_{1/2})\). In the case \(k_0^2=k_1^2=1/2\) we have
\(m_1^-=\lambda_1(-\boldsymbol\delta_{\zeta})=-1/4\) for any \(\zeta\in [0,1]\).
In the case \(k_0^2=k_1^2<1/2\) we have \(m_1^-=\lambda_1(-\boldsymbol\delta_0)\).

%%%%%%%%%%%%%%%%%%%%%%%%%%%%%%%%%%%%%%%%%%%%%%%%%%%%%%%%%%%%%%%%%%%%%%%%%%%%%%%%%%%

\section{The set \protect\(\Gamma_1\protect\) and related topics}\label{par:2}
\subsection{}\label{pt:1.2}
We suppose that all considered functional spaces are real.

By \(W_2^{-1}[0,1]\) denote the Hilbert space that is a completion of \(L_2[0,1]\)
in the norm
\[
	\|y\|_{W_2^{-1}[0,1]}\rightleftharpoons\sup\limits_{\|z\|_{W_2^1[0,1]}=1}
		\int_0^1 yz\,dx.
\]
When \(y\in W_2^{-1}[0,1]\), by \(\int_0^1 yz\,dx\) we sometimes denote the result
\[
	\langle y,z\rangle\rightleftharpoons\lim\limits_{n\to\infty}
		\int_0^1 y_nz\,dx \quad (\text{where } y=\lim_{n\to\infty}y_n,
		\quad y_n\in L_2[0,1])
\]
of applying the linear functional \(y\) to the function \(z\in W_2^1[0,1]\).

For any fixed \(q\in L_1[0,1]\) and \(\lambda\in\mathbb R\) the map taking each
\(y\in W_1^2[0,1]\) satisfying \eqref{eq:eq2} to
\[
	-y''+(q-\lambda)y\in L_1[0,1]
\]
can be extended by continuity to the bounded operator \(T_q(\lambda):
W_2^1[0,1]\to W_2^{-1}[0,1]\). Using integration by part, we get
\begin{equation}\label{eq:kv}
	(\forall y,z\in W_2^1[0,1])\qquad\langle T_q(\lambda)y,z\rangle=
		\int_0^1 \left[y'z'+(q-\lambda)\,yz
		\right]\,dx+k_0^2\,y(0)z(0)+k_1^2\,y(1)z(1).
\end{equation}
Consider the linear operator pencil\footnote{A linear operator pencil \(L\)
is an operator-valued function such that \(L(\lambda)=A+\lambda B\), where
\(\lambda\in\mathbb R\), \(A\) and \(B\) are some operators not depending
on \(\lambda\).} \(T_q:\mathbb R\to \mathcal B(W_2^1[0,1],
W_2^{-1}[0,1])\) that takes any \(\lambda\in\mathbb R\) to the operator
\(T_q(\lambda)\) described by \eqref{eq:kv}. The spectral problem for \(T_q\)
may be considered as a reformulation (or as a generalization in case when
\(q\in W_2^{-1}[0,1]\) is not summable) of boundary value problem \eqref{eq:eq1},
\eqref{eq:eq2}. We can do this due to the following two facts.

\subsubsection{}\label{prop:1.1}
{\itshape For all \(q\in L_1[0,1]\) and \(\lambda\in\mathbb R\) the function
\(y\in W_2^1[0,1]\) belongs to the kernel of the operator \(T_q(\lambda)\), iff
\(y\in W_1^2[0,1]\) and \(y\) is a solution of problem \eqref{eq:eq1},
\eqref{eq:eq2}.
}

\begin{proof}
It directly follows from the definition of the operator \(T_q(\lambda)\) that
for any solution \(y\in W_1^2[0,1]\) of problem \eqref{eq:eq1}, \eqref{eq:eq2}
the equality \(T_q(\lambda)y=0\) holds.

Let us prove the converse. Consider some \(y\in\ker T_q(\lambda)\), and put
\begin{equation}\label{eq:eq3}
	w(x)\rightleftharpoons y'(x)-\int_0^x (q-\lambda)y\,dt.
\end{equation}
For any \(z\in\Wo_2^1[0,1]\), using \eqref{eq:kv}, we have
\begin{equation}\label{eq:eq3l}
	0=\langle T_q(\lambda)y,z\rangle=\int_0^1 wz'\,dx.
\end{equation}
Since the set of the derivatives of all functions \(z\in\Wo_2^1[0,1]\)
is an orthogonal complement in \(L_2[0,1]\) of the set of all constants, from
\eqref{eq:eq3l} it follows that the function \(w\in L_2[0,1]\) is constant.
Combining this with \eqref{eq:eq3}, we get that the function \(y'\) is absolutely
continuous and its generalized derivative equals \((q-\lambda)y\). Now, using
\eqref{eq:kv}, we see that for any \(z\in W_2^1[0,1]\) we get
\[
	0=\langle T_q(\lambda)y,z\rangle=[-y'(0)+k_0^2y(0)]\,z(0)+
		[y'(1)+k_1^2y(1)]\,z(1),
\]
so \(y\) satisfies conditions \eqref{eq:eq2}.
\end{proof}

\subsubsection{}\label{prop:1.2}
{\itshape For any \(q\in W_2^{-1}[0,1]\) the spectrum of the linear operator
pencil \(T_q\) is purely discrete, simple and bounded from below.
}

\begin{proof}
Note that for any \(y\in W_2^1[0,1]\) we have
\[
	\|y^2\|_{W_2^1[0,1]}\leq\sup_{x\in [0,1]}|y(x)|\cdot
		\sqrt{\int_0^1\left[y^2+4(y')^2\right]\,dx}\leq
		2\|y\|_{C[0,1]}\cdot\|y\|_{W_2^1[0,1]},
\]
then, by the embedding theorem, we get
\begin{equation}\label{eq:norms}
	\|y^2\|_{W_2^1[0,1]}\leq C\,\|y\|_{W_2^1[0,1]}^2,
\end{equation}
where \(C\) is some constant.

Since \(C[0,1]\) is densely embedded in \(W_2^{-1}[0,1]\), for any
\(\varepsilon\in (0,1)\) there exists a function \(\tilde q\in C[0,1]\)
such that \(\|\tilde q-q\|_{W_2^{-1}[0,1]}\leq\varepsilon/C\). Using this
and inequality \eqref{eq:norms}, for any \(y\in W_2^1[0,1]\) we get
\begin{equation}\label{eq:1_1_2_2}
	\left|\int_0^1 (\tilde q-q)\,y^2\,dx\right|\leq
	\|\tilde q-q\|_{W_2^{-1}[0,1]}\cdot\|y^2\|_{W_2^1[0,1]}
	\leq\varepsilon\,\|y\|^2_{W_2^1[0,1]}.
\end{equation}
Further, for any \(\kappa>\|\tilde q\|_{C[0,1]}+1\) we have
\(\int_0^1\tilde q y^2\,dx\geq (1-\kappa)\int_0^1 y^2\,dx\). Combining
this with \eqref{eq:kv} and \eqref{eq:1_1_2_2}, we obtain
\begin{flalign}\label{eq:kappa}
	&& IT_q(-\kappa)&\geq 1-\varepsilon,&&
\end{flalign}
where by \(I:W_2^{-1}[0,1]\to W_2^1[0,1]\) we denote an isometry that satisfies
\[
	(\forall y\in W_2^{-1}[0,1])\,(\forall z\in W_2^1[0,1])\qquad
		\langle Iy,z\rangle_{W_2^1[0,1]}=\langle y,z\rangle.
\]
The existence and uniqueness of this isometry follows from the Riesz theorem
about the representation of a functional in a Hilbert space
\cite[\S\S~30,\,99]{RN:1979}.

From estimate \eqref{eq:kappa} it follows \cite[\S~104]{RN:1979} that the operator
\(S\rightleftharpoons IT_q(-\kappa)\) is boundedly invertible.
Taking into account \eqref{eq:kv}, we have \(IT_q(\lambda)\equiv S-(\lambda+
\kappa)J^*J\), where \(J:W_2^1[0,1]\to L_2[0,1]\) is the embedding operator.
So for any \(\lambda\in\mathbb R\) the existence of a bounded inverse of
the operator \(T_q(\lambda)\) is equivalent to the existence of a bounded inverse
of the operator \(1-(\lambda+\kappa)S^{-1/2}J^*JS^{-1/2}\). Since \(J\)
is compact, it follows that the spectrum of \(T_q\) is purely discrete, semi-simple
and bounded from below.

The spectrum of the pencil \(T_q\) is simple since (see \cite{SaSh:2003},
\cite[Propositions~2, 10]{Vl:2009}) for any \(\lambda\in\mathbb R\) the kernel
of the operator \(T_q(\lambda)\) is formed by the first components \(Y_1\)
of the solutions to the boundary value problem
\begin{gather}\label{eq:Matr1}
	\begin{pmatrix}Y_1 \\ Y_2\end{pmatrix}' = \begin{pmatrix}u&1\\
		-u^2&-u\end{pmatrix} \begin{pmatrix}Y_1 \\ Y_2\end{pmatrix},\\
		\label{eq:Matr2}
		Y_2(0)-k_0^2 Y_1(0)=Y_2(1)+[k_1^2+\omega]\,Y_1(1)=0.
\end{gather}
Here \(u\in L_2[0,1]\) and \(\omega\in\mathbb R\) are taken from the representation
\begin{equation}\label{eq:Matr3}
	(\forall y\in W_2^1[0,1])\qquad\int_0^1 (q-\lambda)y\,dx
		=-\int_0^1 uy'\,dx+\omega\,y(1)
\end{equation}
of the potential \(q\in W_2^{-1}[0,1]\).
\end{proof}

\subsection{}
For eigenvalues
\[
	\lambda_1(q)<\lambda_2(q)<\ldots<\lambda_n(q)<\ldots
\]
of the pencil \(T_q\) we have the following propositions.

\subsubsection{}\label{prop:varpr}(see \cite[Proposition~10]{Vl:2009})
{\itshape For any \(n\geq 1\), \(q\in W_2^{-1}[0,1]\) and \(\lambda\in
\mathbb R\) the inequality \(\lambda>\lambda_n(q)\) is equivalent to the existence
of \mbox{\(n\)-di}\-men\-ti\-o\-nal subspace \(\mathfrak N\subset W_2^1[0,1]\)
that satisfies
\[
	\belowdisplayskip=-0.5cm
	(\forall y\in\mathfrak N\setminus\{0\})\qquad
		\langle T_q(\lambda)y,y\rangle<0.
\]
}

\subsubsection{}\label{prop:cont}
{\itshape For any \(n\geq 1\) the function \(\lambda_n:W_2^{-1}[0,1]\to
\mathbb R\) is continuous.
}

\begin{proof}
Consider some \(q\in W_2^{-1}[0,1]\) and \(\varepsilon\in (0,1/2)\). For any
\(y\in W_2^1[0,1]\), \(\lambda\in\mathbb R\) and \(\tilde q\in W_2^{-1}[0,1]\)
such that \(\|\tilde q-q\|_{W_2^{-1}[0,1]}<\varepsilon/C\), where \(C\) is
the same as in \eqref{eq:norms}, we get
\begin{flalign*}
	&& \langle T_{\tilde q}(\lambda)y,y\rangle&\geq\langle T_q(\lambda)y,
		y\rangle-\varepsilon\,\|y\|_{W_2^1[0,1]}^2&&\\
	&& &\geq\langle T_q(\lambda)y,y\rangle-\varepsilon\,
		\|y\|_{W_2^1[0,1]}^2-\varepsilon\cdot\langle T_{2q}(\lambda_1(2q))
		y,y\rangle-\varepsilon k_0^2y^2(0)-\varepsilon k_1^2y^2(1)&&\\
	&& &=(1-2\varepsilon)\cdot\left\langle T_q\left(\dfrac{\lambda+
		\varepsilon\cdot [1-\lambda_1(2q)]}{1-2\varepsilon}\right)y,y
		\right\rangle.
\end{flalign*}
Consequently, from variational principle \ref{prop:varpr} it follows that any
\(\lambda>\lambda_n(\tilde q)\) satisfies
\[
	\dfrac{\lambda+\varepsilon\cdot [1-\lambda_1(2q)]}{1-2\varepsilon}>
		\lambda_n(q).
\]
Since we can choose \(\lambda\) arbitrarily close to \(\lambda_n(\tilde q)\),
we have \(\lambda_n(\tilde q)\geq (1-2\varepsilon)\,\lambda_n(q)-\varepsilon
\cdot[1-\lambda_1(2q)]\). By the same method we get \(\lambda_n(\tilde q)\leq
(1+2\varepsilon)\,\lambda_n(q)+\varepsilon\cdot[1-\lambda_1(2q)]\).
\end{proof}

\subsection{}\label{pt:2.2}
Let \(\Gamma_1\) be the closure in \(W_2^{-1}[0,1]\) of the set \(A_1\)
defined by \eqref{eq:Agamma}. Put by definition
\[
	\Lambda(X)\rightleftharpoons\{\lambda\in\mathbb R\::\:
		(\exists q\in X)\quad\lambda=\lambda_1(q)\},
\]
where \(X\subseteq W_2^{-1}[0,1]\) is some set of generalized functions.
The set \(\Lambda(X)\) is formed by all the possible values of \(\lambda_1(q)\)
for all \(q\in X\). By \(-X\) we, as usually, denote the set
\[
	\{q\in W_2^{-1}[0,1]\::\:(\exists r\in X)\quad q=-r\}.
\]

\subsubsection{}\label{prop:2.2}
{\itshape Suppose \(X\) is a dense subset of\/ \(\Gamma_1\), then the closures
of\/ \(\Lambda(\pm X)\) and\/ \(\Lambda(\pm \Gamma_1)\) coincide.
}

\subsubsection{}\label{prop:2.3.2}
{\itshape The extrema \(m^{\pm}_1\rightleftharpoons\inf\Lambda(\pm A_1)\) and
\(M^{\pm}_1\rightleftharpoons\sup\Lambda(\pm A_1)\), defined in \ref{pt:vrai1.2},
satisfy the equalities \(m^{\pm}_1=\inf\Lambda(\pm\Gamma_1)\) and \(M^{\pm}_1=
\sup\Lambda(\pm\Gamma_1)\).
}

\medskip
Proposition \ref{prop:2.2} immediately follows from \ref{prop:cont}. Proposition
\ref{prop:2.3.2} immediately follows from \ref{prop:2.2}.

\subsubsection{}\label{prop:4.2}
{\itshape The set\/ \(\Gamma_1\) consists of all
non-negative\footnote{The generalized function \(q\in W_2^{-1}[0,1]\) is called
non-negative if for any non-negative function \(y\in W_2^1[0,1]\) the inequality
\(\langle q,y\rangle\geq 0\) holds.} distributions \(q\in W_2^{-1}[0,1]\)
such that \(\int_0^1 q\,dx=1\).
}

\begin{proof}
Since for any \(q\in\Gamma_1\) there exists a sequence of functions from \(A_1\)
such that its limit equals \(q\), it follows that the generalized function \(q\)
is non-negative and satisfies \(\int_0^1 q\,dx=1\).

Let us prove the converse. Suppose \(q\in W_2^{-1}[0,1]\) is a non-negative
generalized function and satisfies \(\int_0^1 q\,dx=1\).
Then (see \cite{SaSh:2003}, \cite[\S\,2.3]{Vl:2009}) there exists a function
\(u\in L_2[0,1]\) such that
\begin{equation}\label{eq:L2}
	(\forall y\in W_2^1[0,1])\qquad\int_0^1 qy\,dx= -\int_0^1 uy'\,dx+y(1).
\end{equation}
Put by definition
\[
	\Pi_{\gamma,\eta,\theta}(x)\rightleftharpoons\left\{\begin{aligned}&
		\tfrac{x-\gamma}{\eta-\gamma}&&\text{for }x\in [\gamma,\eta],\\
		&\tfrac{\theta-x}{\theta-\eta}&&\text{for }x\in [\eta,\theta],\\
		&0&&\text{otherwise}\end{aligned}\right.
\]
for any reals \(\gamma<\eta<\theta\). Suppose \(0<a<b<c<d<1\). Substituting
the functions \(\Pi_{-1,0,a}+\Pi_{0,a,b}\), \(\Pi_{a,b,c}+\Pi_{b,c,d}\) and
\(\Pi_{c,d,1}+\Pi_{d,1,2}\) for \(y\) in \eqref{eq:L2}, we get
\[
	0\leq\dfrac{1}{b-a}\int_a^b u\,dx\leq
		\dfrac{1}{d-c}\int_c^d u\,dx\leq 1.
\]
From these inequalities it follows that the function \(u\in L_2[0,1]\)
is non-decreasing and satisfies \(\vrai{inf}_{x\in [0,1]}u(x)\geq 0\)
and \(\vrai{sup}_{x\in [0,1]}u(x)\leq 1\).

Since there exists a sequence \(\{u_n\}_{n=0}^{\infty}\) of non-decreasing
piecewise linear functions such that \(u_n(0)=0\), \(u_n(1)=1\)
and \(u=\lim\limits_{n\to\infty}{u_n}\), it follows that \(q=\lim\limits_{n\to
\infty}{u'_n}\), where \(u'_n \in A_1\).
\end{proof}

\subsection{}
Consider the function \(F\) implicitely defined by the equation
\begin{equation}\label{eq:implF}
	\lambda_1(F(\mu,\zeta)\boldsymbol{\delta}_{\zeta})=\mu,
\end{equation}
where \(\mu\in\mathbb R\) and \(\zeta\in [0,1]\). The following three propositions
give us some information about this function.

\subsubsection{}\label{prop:F1}
{\itshape For any \(\zeta\in [0,1]\) the function \(F(\cdot,\zeta)\)
is single-valued, strictly increasing, and its domain is the interval \((-\infty,
f^+)\) with some \(f^+>0\).
}

\begin{proof}
For any \(a\in\mathbb R\) there exists \cite[Proposition~11]{Vl:2009} a positive
eigenfunction \(y\in\ker T_{a\boldsymbol{\delta}_{\zeta}}(\mu)\) corresponding
to the eigenvalue \(\mu\rightleftharpoons\lambda_1(a\boldsymbol{\delta}_{\zeta})\),
so for any \(b<a\) we have
\[
	\langle T_{b\boldsymbol{\delta}_{\zeta}}(\mu)y,y\rangle =
		\langle T_{a\boldsymbol{\delta}_{\zeta}}(\mu)y,y\rangle
		+(b-a)\cdot y^2(\zeta)<0.
\]
Using \ref{prop:varpr}, we now get \(\lambda_1(b\boldsymbol{%
\delta}_{\zeta})<\mu\). So the function \(F(\cdot,\zeta)\) is the inverse
of the strictly increasing and, according to \ref{prop:cont}, continuous
map \(a\mapsto \lambda_1(a\boldsymbol\delta_{\zeta})\). Therefore, the function
\(F(\cdot,\zeta)\) is single-valued and strictly increasing.

Further, for any \(a\in\mathbb R\) from the equality
\[
	\langle T_{a\boldsymbol{\delta}_{\zeta}}(a+k_0^2+k_1^2)1,
		1\rangle=a-(a+k_0^2+k_1^2)+k_0^2+k_1^2 =0
\]
and proposition \ref{prop:varpr} it follows that \(\lambda_1(a\boldsymbol{%
\delta}_{\zeta})\leq a+k_0^2+k_1^2\). Therefore the domain of \(F(\cdot,
\zeta)\) is unbounded from below. Also for any \(a>0\) we have \(\lambda_1(a
\boldsymbol\delta_{\zeta})>0\), so the right bound of \(\dom F(\cdot,\zeta)\)
is positive.
\end{proof}

\subsubsection{}\label{prop:F2}
{\itshape The function \(F\) is continuous.
}

\begin{proof}
Consider an arbitrary point \((\mu_0,\zeta_0)\in\dom F\) and suppose \(a^{\pm}\)
satisfy \(a^{-}<F(\mu_0,\zeta_0)<a^{+}\). For any point \((\mu,\zeta)\in\mathbb R
\times [0,1]\) sufficiently close to \((\mu_0,\zeta_0)\) from \ref{prop:F1} and
\ref{prop:cont} we obtain the inequalities \(\lambda_1(a^{-}\boldsymbol{\delta}_{%
\zeta})<\mu<\lambda_1(a^{+} \boldsymbol{\delta}_{\zeta})\). Hence there exists
\(a\in (a^{-},a^+)\) such that \(\mu=\lambda_1(a\boldsymbol{\delta}_{\zeta})\),
so for the point \((\mu,\zeta)\) equation \eqref{eq:implF} has a solution
\(F(\mu,\zeta) = a\).
\end{proof}

\subsubsection{}\label{prop:F3}
{\itshape A point \((\mu,\zeta)\in (0,+\infty)\times [0,1]\) belongs to domain
of the function \(F\) iff the following conditions hold:
\begin{equation}\label{eq:sign}
	\sqrt{\mu}\cdot(\zeta-\alpha_{\mu})\in (-\pi/2,\pi/2),\qquad
	\sqrt{\mu}\cdot(1-\beta_{\mu}-\zeta)\in(-\pi/2,\pi/2),
\end{equation}
where \(\alpha_{\mu}\) and \(\beta_{\mu}\) are defined by \eqref{eq:ab}.
In this case the equality
\begin{equation}\label{eq:F}
	F(\mu,\zeta)=\sqrt{\mu}\cdot\left\{\tan[\sqrt{\mu}\cdot(\zeta-
		\alpha_{\mu})]+\tan[\sqrt{\mu}\cdot(1-\beta_{\mu}-\zeta)]\right\}
\end{equation}
holds.

For any \(\zeta\in [0,1]\) the equality
\begin{equation}\label{eq:F0}
	F(0,\zeta)=-\dfrac{k_0^2}{1+k_0^2\zeta}-\dfrac{k_1^2}{1+k_1^2\,
		(1-\zeta)}
\end{equation}
holds.

For any \(\mu<0\) and \(\zeta\in [0,1]\) the equality
\begin{equation}\label{eq:Fmin}
	F(\mu,\zeta)=-\sqrt{|\mu|}\cdot\{G(\sqrt{|\mu|},k_0^2,\zeta)
		+G(\sqrt{|\mu|},k_1^2,1-\zeta)\},
\end{equation}
where
\[
	G(\nu,\kappa,x)\rightleftharpoons\left\{\begin{aligned}
		&\tanh\left(\nu x+\ln\sqrt{\dfrac{\nu+
		\kappa}{\nu-\kappa}}\right)&&\text{for }\nu>\kappa,\\
		&1&&\text{for }\nu=\kappa,\\
		&\coth\left(\nu x+\ln\sqrt{\dfrac{\kappa+\nu}{\kappa-
		\nu}}\right)&&\text{for }\nu<\kappa,\end{aligned}\right.
\]
holds.
}

\begin{proof}
Consider \(\mu\in\mathbb R\) and \(\zeta\in (0,1)\) such that \((\mu,\zeta)\in
\dom F\). According to \eqref{eq:Matr1}--\eqref{eq:Matr3}, the equality
\(T_{q}(\mu)y=0\), where \(q\rightleftharpoons F(\mu,\zeta)
\boldsymbol\delta_{\zeta}\), is equivalent to the boundary problem
\begin{gather}\label{eq:delt1}
	-y''=\mu y\quad\text{at }(0,\zeta)\cup (\zeta,1),\\
	\label{eq:delt2} y'(\zeta+0)-y'(\zeta-0)=F(\mu,\zeta)y(\zeta),\\
	\label{eq:delt3} y'(0)-k_0^2y(0)=y'(1)+k_1^2y(1)=0.
\end{gather}
From \cite[Proposition~11]{Vl:2009} and \eqref{eq:implF} it follows that
any non-trivial solution to problem \eqref{eq:delt1}--\eqref{eq:delt3}
is constant-sign.

In the case \(\mu>0\) any solution to problem \eqref{eq:delt1}, \eqref{eq:delt3}
has the form
\begin{equation}\label{eq:2.19}
	y(x)=\left\{\begin{aligned}&A\cdot\cos[\sqrt{\mu}\cdot(1-\beta_{\mu}-\zeta)]
		\cdot\cos[\sqrt{\mu}\cdot(x-\alpha_{\mu})]&&\text{for }x<\zeta,\\
		&A\cdot\cos[\sqrt{\mu}\cdot(1-\beta_{\mu}-x)]\cdot\cos[\sqrt{\mu}
		\cdot(\zeta-\alpha_{\mu})]&&\text{for }x>\zeta,\end{aligned}\right.
\end{equation}
where \(A\) is some constant. This function is constant-sign iff conditions
\eqref{eq:sign} hold. Using \eqref{eq:delt2}, we now get \eqref{eq:F}. The values
\(\zeta\in\{0,1\}\) are finally included in the consideration using propositions
\ref{prop:F2} and \ref{prop:cont}.

The cases \(\mu=0\) and \(\mu<0\) are considered on the base of
\eqref{eq:delt1}--\eqref{eq:delt3} by analogous way using the solution
\begin{equation}\label{eq:2.20}
	y(x)=\left\{\begin{aligned}&A\cdot [1+k_1^2(1-\zeta)]\cdot[1+k_0^2x]
		&&\text{for }x<\zeta,\\ &A\cdot [1+k_1^2(1-x)]\cdot[1+k_0^2\zeta]
		&&\text{for }x>\zeta\end{aligned}\right.
\end{equation}
in the case \(\mu=0\), and the solution
\begin{equation}\label{eq:2.21}
	y(x)=\left\{\begin{aligned}&A\cdot g(\sqrt{|\mu|},k_1^2,1-\zeta)\cdot
		g(\sqrt{|\mu|},k_0^2,x)&&\text{for }x<\zeta,\\
		&A\cdot g(\sqrt{|\mu|},k_1^2,1-x)\cdot g(\sqrt{|\mu|},k_0^2,
		\zeta)&&\text{for }x>\zeta,\end{aligned}\right.
\end{equation}
where
\[
	g(\nu,\kappa,x)\rightleftharpoons\left\{\begin{aligned}
		&\cosh\left(\nu x+\ln\sqrt{\dfrac{\nu+
		\kappa}{\nu-\kappa}}\right)&&\text{for }\nu>\kappa,\\
		&e^{\nu x}&&\text{for }\nu=\kappa,\\
		&\sinh\left(\nu x+\ln\sqrt{\dfrac{\kappa+\nu}{\kappa-
		\nu}}\right)&&\text{for }\nu<\kappa,\end{aligned}\right.
\]
in the case \(\mu<0\).
\end{proof}

%%%%%%%%%%%%%%%%%%%%%%%%%%%%%%%%%%%%%%%%%%%%%%%%%%%%%%%%%%%%%%%%%%%%%%%%%%%%%

\section{Proofs of the main results}\label{par:3}
\subsection{}
In this section we prove theorems \ref{prop:M1+}--\ref{prop:m1-}.
We use the notation
\begin{align*}
	\Omega^+(y)&\rightleftharpoons\textstyle\left\{x\in[0,1]\::\:
		y(x)=\sup_{t\in[0,1]}y(t)\right\},\\
	\Omega^-(y)&\rightleftharpoons\textstyle\left\{x\in[0,1]\::\:
		 y(x)=\inf_{t\in[0,1]}y(t)\right\},
\end{align*}
where \(y\in W_2^1[0,1]\) is an arbitrary positive function. Also we take
into account proposition \ref{prop:2.3.2}.

\subsection{Proof of theorem \ref{prop:M1+}.}
Consider some potential \(q^*\in\Gamma_1\), and some positive eigenfunction
\(y\in\ker T_{q^*}(\lambda_1(q^*))\). Suppose that the support of the generalized
function \(q^*\) is a subset of \(\Omega^+(y)\). Then for any \(q\in\Gamma_1\)
we, using \ref{prop:4.2}, have
\begin{flalign*}
	&&0&=\langle T_{q^*}(\lambda_1(q^*))y,y\rangle\\
	&& &=\int_0^1\left[(y')^2- \lambda_1(q^*)\,y^2\right]\,dx+
		\sup\limits_{x\in [0,1]}y^2(x)+k_0^2\,y^2(0)+k_1^2\,y^2(1)\\
	&& &\geq\int_0^1\left[(y')^2+(q-\lambda_1(q^*))\,y^2\right]\,dx+
		k_0^2\,y^2(0)+k_1^2\,y^2(1),&&
\end{flalign*}
hence \(\langle T_q(\lambda_1(q^*))y,y\rangle\leq 0\). It follows that
\(\lambda_1(q)\leq\lambda_1(q^*)\), therefore \(\lambda_1(q^*)=M_1^+\).
Thus we proved that \(M_1^+\) is attained at any potential \(q^*\) such that
\(\supp q^*\subseteq\Omega^+(y)\).

Suppose that \(\Omega^+(y) = [\tau_0,\tau_1]\), where \(\tau_0 \ne \tau_1\).
Also suppose that the potential \(q^*\) is summable and has the form
\[
	q^*(x)= \left\{\begin{aligned}&\mu& &\text{for }x\in[\tau_0,\,\tau_1],\\
		&0& &\text{otherwise},\end{aligned}\right.
\]
where \(\mu\) is some positive constant. Since \(y''(x)=0\) for all \(x\in
(\tau_0,\tau_1)\), it follows that \(\mu=\lambda_1(q^*)\).
Therefore the eigenfunction \(y\) has the form
\[
	y(x)=\left\{\begin{aligned}
		&A\cdot\cos[\sqrt{\mu}\cdot(x-\alpha_{\mu})]&
		&\text{for }x<\tau_0,\\
		&B&&\text{for }x\in[\tau_0,\,\tau_1],\\
		&C\cdot\cos[\sqrt{\mu}\cdot(1-\beta_{\mu}-x)]&
		&\text{for }x>\tau_1,\end{aligned}\right.
\]
where \(A\), \(B\) and \(C\) are some positive constants, and \(\alpha_{\mu}\),
\(\beta_{\mu}\) are defined by \eqref{eq:ab}. From the continuity of \(y'\)
it follows that \(\tau_0=\alpha_{\mu}\) and \(\tau_1=1-\beta_{\mu}\), hence
\(A=B=C\). Finally, from the condition \(\int_0^1 q^*\,dx=1\) we have
the equation \eqref{eq:M1}.

To conclude the proof, it remains to note that equation \eqref{eq:M1} has
a unique solution, because \(\alpha_{\mu}\) and \(\beta_{\mu}\), considered
as functions of \(\mu>0\), are non-negative, continuous, non-increasing
and tend to zero as \(\mu\to+\infty\).

\subsection{Proof of theorem \ref{prop:M1-}.}
Consider some potential \(q^*\in-\Gamma_1\), and some positive eigenfunction
\(y\in\ker T_{q^*}(\lambda_1(q^*))\). Suppose that \(\supp q^*\subseteq
\Omega^-(y)\). Then for any \(q\in-\Gamma_1\) we, using \ref{prop:4.2}, have
\begin{flalign*}
	&&0&=\langle T_{q^*}(\lambda_1(q^*))y,y\rangle\\
	&& &=\int_0^1\left[(y')^2- \lambda_1(q^*)\,y^2\right]\,dx-
		\inf\limits_{x\in [0,1]}y^2(x)+k_0^2\,y^2(0)+k_1^2\,y^2(1)\\
	&& &\geq\int_0^1\left[(y')^2+(q-\lambda_1(q^*))\,y^2\right]\,dx+
		k_0^2\,y^2(0)+k_1^2\,y^2(1),&&
\end{flalign*}
hence \(\langle T_q(\lambda_1(q^*))y,y\rangle\leq 0\). It follows that
\(\lambda_1(q)\leq\lambda_1(q^*)\), therefore \(\lambda_1(q^*)=M_1^-\).
Thus we proved that \(M_1^-\) is attained at any potential \(q^*\) such that
\(\supp q^*\subseteq\Omega^-(y)\).

Suppose \(k_0^2+k_1^2\leq 1\). Consider the generalized function
\(q^*\rightleftharpoons -k_0^2\boldsymbol{\delta}_0-k_1^2\boldsymbol{\delta}_1-
(1-k_0^2-k_1^2)\), which in this case belongs to \(-\Gamma_1\). Using
\eqref{eq:kv}, we get that the first eigenfunction of the pencil \(T_{q^*}\)
is \(y\equiv\mathrm{const}\), so \(\supp{q^*}\subseteq\Omega^-(y)\). It follows
that \(M_1^-\) is attained at the potential \(q^*\) and is equal to
the corresponding first eigenvalue \(\lambda_1(q^*)=k_0^2+k_1^2-1\).

Suppose
\begin{gather}\label{eq:first}
	k_0^2+k_1^2\geq 1,\\ \label{eq:second}
	k_1^2-k_0^2\leq 1.
\end{gather}
Consider the generalized function \(q^*\rightleftharpoons-(1+k_0^2-k_1^2)
\boldsymbol{\delta}_0/2-(1-k_0^2+k_1^2)\boldsymbol{\delta}_1/2\), which, due to
\eqref{eq:second}, belongs to \(-\Gamma_1\). For such \(q^*\) the equation
\(T_{q^*}(\lambda)y=0\) is equivalent to problem \eqref{eq:propM1-1},
\eqref{eq:propM1-2}. The first eigenvalue \(\lambda_1(q^*)\), due to
\eqref{eq:first} and \eqref{eq:propM1-2}, is non-negative and the corresponding
eigenfunction is
\begin{equation}\label{eq:100500}
	y(x)\equiv\cos[\sqrt{\lambda_1(q^*)}\cdot (x-\zeta)],
\end{equation}
where \(\zeta=1/2\). Hence \(\supp q^*\subseteq\Omega^-(y)\). It follows that
\(M_1^-\) is attained at the potential \(q^*\) and is equal to the corresponding
first eigenvalue \(\lambda_1(q^*)\).

Suppose \(k_1^2-k_0^2\geq 1\). Consider the generalized function
\(q^*\rightleftharpoons-\boldsymbol{\delta}_1\in-\Gamma_1\). For such \(q^*\)
the equation \(T_{q^*}(\lambda)y=0\) is equivalent to problem \eqref{eq:propM1-1},
\eqref{eq:propM1-3}. The corresponding first eigenfunction is defined by
\eqref{eq:100500}, where \(\zeta\in [0,1/2]\), since \(k_1^2-1\geq k_0^2\).
Hence \(\supp q^*\subseteq\Omega^{-}(y)\). It follows that \(M_1^-\) is attained
at the potential \(q^*\) and is equal to the corresponding first eigenvalue
\(\lambda_1(q^*)\).

\subsection{Proof of theorem \ref{prop:m1+}.}
Consider some potential \(q\in\Gamma_1\), and some positive eigenfunction \(y\in
\ker T_{q}(\lambda_1(q))\). Then for any \(\lambda>\lambda_1(q)\), according to
\ref{prop:4.2}, we have
\begin{flalign*}
	&& 0&>\int_0^1 \left[(y')^2+(q-\lambda)\,y^2\right]\,dx
		+k_0^2\,y^2(0)+k_1^2\,y^2(1)\\
	&& &\geq\int_0^1 \left[(y')^2-\lambda y^2\,\right]\,dx
		+\inf\limits_{x\in[0,1]} y^2(x)+k_0^2\,y^2(0)
		+k_1^2\,y^2(1).&
\end{flalign*}
It follows that there exists \(\zeta\in [0,1]\) such that
\[
	\int_0^1 \left[(y')^2+(\boldsymbol{\delta}_{\zeta}-\lambda)\,
		y^2\right]\,dx+k_0^2\,y^2(0)+k_1^2\,y^2(1)<0.
\]
So for any \(\lambda>m^+_1\) there exists \(\zeta\in [0,1]\) such that
\(\lambda_1(\boldsymbol{\delta}_{\zeta})<\lambda\). Hence, using \ref{prop:4.2},
we get \(m^+_1=\inf_{x\in [0,1]}\lambda_1(\boldsymbol{\delta}_x)\). This equality
is equivalent, according to \ref{prop:F1}, to the following fact: \(F(m^+_1,x)\)
is defined for all \(x\in [0,1]\) and satisfies \(\sup_{x\in [0,1]}F(m^+_1,x)=1\).

Since \(m_1^+>0\), from \ref{prop:F3} it follows that if \(\mu=m^+_1\), then for
any \(\zeta\in [0,1]\) conditions \eqref{eq:sign} hold. According to \eqref{eq:F},
\eqref{eq:sign} and
\begin{equation}\label{eq:diff}
	\dfrac{\partial F(\mu,\zeta)}{\partial\zeta}\equiv\mu\cdot\dfrac{%
		\cos^2[\sqrt{\mu}\cdot(1-\beta_{\mu}-\zeta)]-\cos^2[\sqrt{\mu}
		\cdot(\zeta-\alpha_{\mu})]}{\cos^2[\sqrt{\mu}\cdot(\zeta-
		\alpha_{\mu})]\cdot\cos^2[\sqrt{\mu}\cdot(1-\beta_{\mu}-\zeta)]},
\end{equation}
it follows that the function \(F(\mu,\cdot)\) can have at some point \(\zeta\in
(0,1)\) a local extremum satisfying \(F(\mu,\zeta)>0\) only if \(\zeta=(1-
\beta_{\mu}+\alpha_{\mu})/2\), \(\zeta>\alpha_{\mu}\) and \(\zeta<1-\beta_{\mu}\).
But this conditions imply, according to \eqref{eq:diff}, that such \(\zeta\)
must be a point of strict local minimum of the function \(F(\mu,\cdot)\). Therefore
\(F(\mu,\cdot)\) can't have a supremum in \((0,1)\), so we get \(m^+_1=\inf\{%
\lambda_1(\boldsymbol{\delta}_0),\lambda_1(\boldsymbol{\delta}_1)\}\). Note that
for the potential \(q^*\rightleftharpoons\boldsymbol\delta_i\), where \(i\in
\{0,1\}\), the equation \(T_{q^*}(\lambda)y=0\) is equivalent to the problem
\begin{gather*}
	-y''=\lambda y,\\ y'(0)-[k_0^2+(1-i)]\,y(0)=y'(1)+[k_1^2+i]\,y(1)=0.
\end{gather*}
Therefore we have
\[
	\dfrac{\lambda_1(\boldsymbol{\delta}_i)-k_0^2k_1^2-k_{1-i}^2}{%
		k_0^2+k_1^2+1}=\sqrt{\lambda_1(\boldsymbol{\delta}_i)}
		\cot\sqrt{\lambda_1(\boldsymbol{\delta}_i)},
\]
so \(m^+_1=\lambda_1(\boldsymbol{\delta}_1)\).

\subsection{Proof of theorem \ref{prop:m1-}.}
Consider some potential \(q\in-\Gamma_1\), and some positive eigenfunction
\(y\in\ker T_{q}(\lambda_1(q))\). Then for any \(\lambda>\lambda_1(q)\), according
to \ref{prop:4.2}, we have
\begin{flalign*}
	&& 0&>\int_0^1 \left[(y')^2+(q-\lambda)\,y^2\right]\,dx
		+k_0^2\,y^2(0)+k_1^2\,y^2(1)\\
	&& &\geq\int_0^1 \left[(y')^2-\lambda y^2\,\right]\,dx
		-\sup\limits_{x\in[0,1]} y^2(x)+k_0^2\,y^2(0)
		+k_1^2\,y^2(1).&&
\end{flalign*}
It follows that there exists \(\zeta\in [0,1]\) such that
\[
	\int_0^1 \left[(y')^2+(-\boldsymbol{\delta}_{\zeta}-\lambda)\,
		y^2\right]\,dx+k_0^2\,y^2(0)+k_1^2\,y^2(1)<0.
\]
So for any \(\lambda>m^-_1\) there exists \(\zeta\in [0,1]\) such that
\(\lambda_1(-\boldsymbol{\delta}_{\zeta})<\lambda\). Hence, using \ref{prop:4.2},
we get \(m^-_1=\inf_{x\in [0,1]}\lambda_1(-\boldsymbol{\delta}_x)\). This equality
is equivalent, according to \ref{prop:F1}, to the following fact: \(F(m^-_1,x)\)
is defined for all \(x\in [0,1]\) and satisfies \(\sup_{x\in [0,1]}F(m^-_1,x)=-1\).

For any fixed value \(\mu\in\mathbb R\) we consider the conditions
\begin{gather}\label{eq:Fstat1}
	F(\mu,\zeta)<0,\\ \label{eq:Fstat2}
	\partial F(\mu,\zeta)/\partial\zeta=0.
\end{gather}
It is clear that some point \(\zeta\in (0,1)\) can satisfy the equalities
\(F(\mu,\zeta)=\sup_{x\in [0,1]} F(\mu,x)=-1\) only if \eqref{eq:Fstat1}
and \eqref{eq:Fstat2} hold.

Suppose \(\mu>0\). Then, according to \eqref{eq:diff}, \eqref{eq:F},
\eqref{eq:2.19} and \eqref{eq:delt2}, for any point \(\zeta\in (0,1)\)
satisfying \eqref{eq:Fstat1} condition \eqref{eq:Fstat2} holds iff the problem
\begin{gather}\label{eq:grmm3}
	-y''=\mu y\quad\text{at } (0,\zeta)\cup(\zeta,1),\\ \label{eq:grmm4}
	\hbox to 0.8\textwidth{$y'(0)-k_0^2y(0)=2y'(\zeta-0)+
		F(\mu,\zeta)y(\zeta)=$\hfill}\\ \notag
	\hbox to 0.8\textwidth{\hfill$=2y'(\zeta+0)-F(\mu,\zeta)y(\zeta)=
		y'(1)+k_1^2y(1)=0$}
\end{gather}
has a continuous positive solution. Besides, for any point \(\zeta\in (0,1)\)
satisfying \eqref{eq:Fstat1} and \eqref{eq:Fstat2} we have
\begin{equation}\label{eq:39}
	\alpha_{\mu}>\zeta>1-\beta_{\mu}.
\end{equation}
Therefore, according to \eqref{eq:diff} and \eqref{eq:sign}, this stationary point
\(\zeta\) is a strict maximum of \(F(\mu,\cdot)\). Since for any \(x\in [0,1]\),
using \eqref{eq:39}, we get
\[
	-\pi/2<-\sqrt{\mu}\alpha_{\mu}\leq\sqrt{\mu}\cdot(x-\alpha_{\mu})
	<\sqrt{\mu}\cdot(x-1+\beta_{\mu})\leq\sqrt{\mu}\beta_{\mu}<\pi/2,
\]
it follows from proposition \ref{prop:F3} that the function \(F(\mu,\cdot)\)
is defined everywhere on \([0,1]\).

Suppose \(\mu=0\). Let us use the same method as in previous case, changing
\eqref{eq:F} to \eqref{eq:F0}, and \eqref{eq:2.19} to \eqref{eq:2.20}. Then we get
that for any point \(\zeta\in (0,1)\) satisfying \eqref{eq:Fstat1} condition
\eqref{eq:Fstat2} holds iff problem \eqref{eq:grmm3}, \eqref{eq:grmm4} has
a continuous positive solution. Using \eqref{eq:F0} we also get that the second
derivative of \(F(0,\cdot)\) is negative. Hence any stationary point \(\zeta\in
(0,1)\) is a strict maximum of \(F(0,\cdot)\).

Suppose \(\mu\in(-k_0^4,0)\). Then, using \eqref{eq:Fmin}, we get
\[
	\dfrac{\partial F(\mu,\zeta)}{\partial\zeta}\equiv-\mu\left\{
		\sinh^{-2}\left(\sqrt{|\mu|}\zeta+\alpha_{\mu}\right)
		-\sinh^{-2}\left(\sqrt{|\mu|}(1-\zeta)+\beta_{\mu}\right)\right\},
\]
where
\[
	\alpha_{\mu}\rightleftharpoons\dfrac12\,\ln\dfrac{k_0^2+\sqrt{|\mu|}}{%
		k_0^2-\sqrt{|\mu|}},\qquad
	\beta_{\mu}\rightleftharpoons\dfrac12\,\ln\dfrac{k_1^2+\sqrt{|\mu|}}{%
		k_1^2-\sqrt{|\mu|}}.
\]
Therefore, according to \eqref{eq:2.21}, for any point \(\zeta\in (0,1)\)
satisfying \eqref{eq:Fstat1} condition \eqref{eq:Fstat2} holds iff problem
\eqref{eq:grmm3}, \eqref{eq:grmm4} has a continuous positive solution. Since
\(\partial^2 F(\mu,\zeta)/\partial\zeta^2<0\), it follows that any stationary
point \(\zeta\in (0,1)\) is a strict maximum of \(F(\mu,\cdot)\).

Suppose \(0>\mu=-k_0^4=-k_1^4\). Then the function \(F(\mu,\cdot)\) is a negative
constant, and for any point \(\zeta\in (0,1)\) problem \eqref{eq:grmm3},
\eqref{eq:grmm4} has a continuous positive solution.

Suppose \(\mu\in [-k_1^4,-k_0^4]\), also \(\mu<0\) and \(k_1>k_0\). Then from
\eqref{eq:Fmin} and \eqref{eq:2.21} it follows that \({\partial F(\mu,\zeta)}/
{\partial\zeta}<0\), and problem \eqref{eq:grmm3}, \eqref{eq:grmm4} has no positive
solutions for any \(\zeta\in (0,1)\).

Suppose \(\mu<-k_1^4\). Then, using \eqref{eq:Fmin}, we get
\[
	\dfrac{\partial F(\mu,\zeta)}{\partial\zeta}\equiv\mu\left\{
		\cosh^{-2}\left(\sqrt{|\mu|}\zeta+\alpha_{\mu}\right)
		-\cosh^{-2}\left(\sqrt{|\mu|}(1-\zeta)+\beta_{\mu}\right)\right\},
\]
where
\[
	\alpha_{\mu}\rightleftharpoons\dfrac12\,\ln\dfrac{\sqrt{|\mu|}+k_0^2}{%
		\sqrt{|\mu|}-k_0^2},\qquad
	\beta_{\mu}\rightleftharpoons\dfrac12\,\ln\dfrac{\sqrt{|\mu|}+k_1^2}{%
		\sqrt{|\mu|}-k_1^2}.
\]
Therefore, according to \eqref{eq:2.21}, for any point \(\zeta\in (0,1)\)
satisfying \eqref{eq:Fstat1} condition \eqref{eq:Fstat2} holds iff problem
\eqref{eq:grmm3}, \eqref{eq:grmm4} has a continuous positive solution. Since
\(\partial^2 F(\mu,\zeta)/\partial\zeta^2>0\), it follows that any stationary
point \(\zeta\in (0,1)\) is a strict minimum of \(F(\mu,\cdot)\).

From proposition \ref{prop:F1} we also get that for any \(\mu\leq 0\) the function
\(F(\mu,\cdot)\) is defined everywhere on \([0,1]\).

Combining all this, we obtain the following: the existence of a continuous positive
solution to problem \eqref{eq:grmm1}, \eqref{eq:grmm2} for some \(\mu\geq-k_0^4\)
and \(\zeta\in (0,1)\) implies that \(F(\mu,\zeta)=-1\), the function \(F(\mu,
\cdot)\) is defined everywhere on \([0,1]\), and \(\sup_{x\in [0,1]}F(\mu,x)\leq
-1\). Therefore \(m_1^-=\lambda_1(-\boldsymbol\delta_{\zeta})\). In converse,
if for any \(\mu\geq-k_0^4\) and \(\zeta\in (0,1)\) the positive solution of
\eqref{eq:grmm1}, \eqref{eq:grmm2} doesn't exist, we get \(m^{-}_1=
\inf\{\lambda_1(-\boldsymbol{\delta}_0),\lambda_1(-\boldsymbol{\delta}_1)\}\).
From the equation
\[
	\lambda_1(-\boldsymbol{\delta}_i)-k_0^2k_1^2+k_{1-i}^2
		=(k_0^2+k_1^2-1)\cdot\psi(\lambda_1(-\boldsymbol{\delta}_i)),
\]
where \(i\in\{0,1\}\) and
\[
	\psi(x)\rightleftharpoons\left\{\begin{aligned}
		&\sqrt{x}\cot\sqrt{x}&&\text{for }x>0,\\
		&1&&\text{for }x=0,\\
		&\sqrt{|x|}\coth\sqrt{|x|}&&\text{for }x<0,\end{aligned}\right.
\]
we obtain that \(\inf\{\lambda_1(-\boldsymbol{\delta}_0),
\lambda_1(-\boldsymbol{\delta}_1)\}=\lambda_1(-\boldsymbol{\delta}_0)\).

\subsection{}\label{pt:3.5}
Now we get some conditions for the existence of a continuous positive solution
to problem \eqref{eq:grmm1}, \eqref{eq:grmm2} considered in theorem \ref{prop:m1-}.

Suppose \(\mu_0(\zeta)\), where \(\zeta\in (0,1]\), is the minimal eigenvalue
of the problem
\begin{gather}\label{eq:grmm5}
	-y''=\lambda y,\\ \label{eq:grmm6}
	y'(0)-k_0^2y(0)=2y'(\zeta)-y(\zeta)=0,
\end{gather}
and suppose \(\mu_1(\zeta)\), where \(\zeta\in [0,1)\), is the minimal eigenvalue
of problem \eqref{eq:grmm5} and
\[
	2y'(\zeta)+y(\zeta)=y'(1)+k_1^2y(1)=0.
\]
It is clear that for some \(\mu\in\mathbb R\) and \(\zeta\in (0,1)\) a continuous
positive solution to \eqref{eq:grmm1}, \eqref{eq:grmm2} exists iff the equalities
\(\mu_0(\zeta)=\mu_1(\zeta)=\mu\) hold.

\subsubsection{}\label{prop:3.2.1}
{\itshape If \(k_0^2=1/2\), then \(\mu_0(\zeta)\equiv-1/4\).

If \(k_0^2>1/2\), then the function \(\mu_0\) strictly decreases and satisfies
\(\lim\limits_{\zeta\to 0}\mu_0(\zeta)=+\infty\) and \(\mu_0(1)>-1/4\).

If \(k_0^2<1/2\), then for any \(\zeta\in (0,1]\) the inequality
\(\mu_0(\zeta)<-1/4\) holds.
}

\begin{proof}
Suppose \(k_0^2=1/2\). Then for any \(\zeta\in (0,1]\) problem
\eqref{eq:grmm5}, \eqref{eq:grmm6} has the positive eigenfunction \(y(x)\equiv
e^{x/2}\) corresponding to the eigenvalue \(-1/4\).

Suppose \(k_0^2>1/2\). Since the eigenvalues of problem \eqref{eq:grmm5},
\eqref{eq:grmm6} increase by \(k_0^2\), it follows that \(\mu_0(\zeta)>-1/4\).
Then let \(y_0\in W_2^1[0,\zeta]\) be an eigenfunction of problem
\eqref{eq:grmm5}, \eqref{eq:grmm6} corresponding to the eigenvalue
\(\mu_0(\zeta)\). Continuing the function \(y_0\) for any \(\theta\in (\zeta,1]\)
to the interval \((\zeta,\theta]\) in the form \(y(x)\rightleftharpoons y_0(\zeta)
e^{(x-\zeta)/2}\), for the obtained function \(y\in W_2^1[0,\theta]\) we get
\[
	\int_0^{\theta} \left[(y')^2-\mu_0(\zeta)\,y^2\right]\,dx+
		k_0^2\,y^2(0)-\dfrac{y^2(\theta)}{2}=[-1/4-\mu_0(\zeta)]\cdot
		[e^{\theta-\zeta}-1]\cdot y^2(\zeta)<0,
\]
hence \(\mu_0(\theta)<\mu_0(\zeta)\). Finally, for \(\zeta\to 0\) we have uniform
by \(y\in W_2^1[0,\zeta]\) asymptotic estimate
\begin{flalign*}
	&& \int_0^{\zeta} (y')^2\,dx+k_0^2\,y^2(0)-\dfrac{y^2(\zeta)}{2}&=
	\left[\int_0^{\zeta}\dfrac{(y')^2}{2}\,dx+(k_0^2-1/2)\,y^2(0)\right]
		{+}\left[\int_0^{\zeta}\dfrac{(y')^2}{2}\,dx+\dfrac{y^2(0)-
		y^2(\zeta)}{2}\right]\\
	&& &\geq\dfrac{k_0^2-1/2+o(1)}{\zeta}\int_0^{\zeta} y^2\,dx-
		\dfrac{1}{2}\int_0^{\zeta} y^2\,dx,
\end{flalign*}
therefore \(\mu_0(\zeta)\geq [k_0^2-1/2+o(1)]\cdot\zeta^{-1}\).

The inequality \(\mu_0(\zeta)<-1/4\) for the case \(k_0^2<1/2\) is proved
likewise  the inequality  \(\mu_0(\zeta)>-1/4\) for the case \(k_0^2>1/2\).
\end{proof}

\subsubsection{}\label{prop:3.2.2}
{\itshape If \(k_1^2=1/2\), then \(\mu_1(\zeta)\equiv-1/4\).

If \(k_1^2>1/2\), then the function \(\mu_1\) strictly increases and satisfies
\(\lim\limits_{\zeta\to 1}\mu_1(\zeta)=+\infty\) and \(\mu_1(0)>-1/4\).

If \(k_1^2<1/2\), then for any \(\zeta\in [0,1)\) the inequality
\(\mu_1(\zeta)<-1/4\) holds.
}

\medskip
The proposition \ref{prop:3.2.2} is proved likewise \ref{prop:3.2.1}.

Combining \ref{prop:3.2.1} and \ref{prop:3.2.2}, we get the last proposition:

\subsubsection{}\label{prop:3.2.3}
{\itshape Problem \eqref{eq:grmm1}, \eqref{eq:grmm2} has a continuous positive
solution for some \(\mu\geq -k_0^4\) and \(\zeta\in (0,1)\) iff
one of these two conditions holds: \(k_0^2>1/2\) or \(k_0^2=k_1^2=1/2\).
}


\begin{thebibliography}{99}
\bibitem{EK:1996} Yu.\,V.~Egorov, V.\,A.~Kondratiev. \emph{On Spectral theory
of elliptic operators}// Operator theory: Advances and Applications, V.~89.
Birk\-hou\-ser, 1996.
\bibitem{VS:2003} V.\,A.~Vinokurov, V.\,A.~Sadovnichii. \emph{On the range
of variation of an eigenvalue when potential is varied}// Dokl. Math. --- 2003. ---
V.~68\,(2). --- P.~247--252.
\bibitem{Ez:2005} S.\,S.~Ezhak. \emph{On the estimates for the minimum eigenvalue
of the Sturm--Liouville problem with integral condition}// Journ. Math. Sci. ---
2007. --- V.~145\,(5). --- P.~5205--5218.
\bibitem{Mur:2005} O.\,V.~Muryshkina. \emph{Estimates for the Minimal Eigenvalue
of the Sturm--Liouville Problem with Nonsymmetric Boundary Conditions}//
Diff. Equat. --- 2001. --- V.~37\,(6). --- P.~899--900.
\bibitem{Kar:2011}  E.~Karulina. \emph{Some estimates for the minimal eigenvalue
of the Sturm--Liouville problem with third-type boundary conditions}//
Mathematica Bohemica. --- 2011. --- V.~136. --- No.~4. --- P.~377--384.
\bibitem{SaSh:2003} A.\,M.~Savchuk, A.\,A.~Shkalikov. \emph{Sturm--Liouville
operators with distribution potentials}// Trans. Mosc. Math. Soc. --- 2003. ---
V.~64. --- P.~143--192.
\bibitem{RN:1979} F.~Riesz, B.~Sz.-Nagy. \emph{Le\c cons d'analyse fonctionelle}.
Budapest: Akad\'emiai Kiad\'o, 1968.
\bibitem{Vl:2009} A.\,A.~Vladimirov. \emph{On the oscillation theory
of the Sturm--Liouville problem with singular coefficients}// Comput. Math.
and Mathem. Physics --- 2009. --- V.~49\,(9). --- P.~1535-1546.
\end{thebibliography}
\end{document}